\documentclass[twoside]{amsart}
\usepackage{url,lastpage}
\usepackage{ifpdf}
\ifpdf
 \usepackage[hyperindex,pagebackref]{hyperref}%
\else
 \expandafter\ifx\csname dvipdfm\endcsname\relax
 \usepackage[hypertex,hyperindex,pagebackref]{hyperref}
 \else
 \usepackage[dvipdfm,hyperindex,pagebackref]{hyperref}
 \fi
\fi
\allowdisplaybreaks[4]
\theoremstyle{plain}
\newtheorem{theorem}{Theorem}[section]
\theoremstyle{remark}
\newtheorem{remark}{Remark}[section]
\DeclareMathOperator{\td}{d}
\newcommand{\bell}{\textup{B}}
\numberwithin{equation}{section}

\begin{document}

\title[An explicit formula for Bernoulli numbers]
{Alternative proofs of a formula for Bernoulli numbers in terms of Stirling numbers}

\author[B.-N. Guo]{Bai-Ni Guo}
\address[Guo]{School of Mathematics and Informatics, Henan Polytechnic University, Jiaozuo City, Henan Province, 454010, China}
\email{\href{mailto: B.-N. Guo <bai.ni.guo@gmail.com>}{bai.ni.guo@gmail.com}, \href{mailto: B.-N. Guo <bai.ni.guo@hotmail.com>}{bai.ni.guo@hotmail.com}}

\author[F. Qi]{Feng Qi}
\address[Qi]{School of Mathematics and Informatics, Henan Polytechnic University, Jiaozuo City, Henan Province, 454010, China}
\email{\href{mailto: F. Qi <qifeng618@gmail.com>}{qifeng618@gmail.com}, \href{mailto: F. Qi <qifeng618@hotmail.com>}{qifeng618@hotmail.com}, \href{mailto: F. Qi <qifeng618@qq.com>}{qifeng618@qq.com}}
\urladdr{\url{http://qifeng618.wordpress.com}}

\begin{abstract}
In the paper, the authors provide four alternative proofs of an explicit formula for computing Bernoulli numbers in terms of Stirling numbers of the second kind.
\end{abstract}

\keywords{alternative proof; explicit formula; Bernoulli numbers; Stirling numbers of the second kind; Fa\`a di Bruno formula; Bell polynomial}

\subjclass[2010]{Primary 11B68, Secondary 11B73}

\thanks{Please cite this article as ``Bai-Ni Guo and Feng Qi, \textit{Alternative proofs of a formula for Bernoulli numbers in terms of Stirling numbers}, Analysis---International mathematical journal of analysis and its applications \textbf{34} (2014), no.~2, 187\nobreakdash--193; Available online at \url{http://dx.doi.org/10.1515/anly-2012-1238}.''}

\maketitle

\section{Introduction}

It is well known that Bernoulli numbers $B_{k}$ for $k\ge0$ may be generated by
\begin{equation}\label{Bernumber-dfn}
\frac{x}{e^x-1}=\sum_{k=0}^\infty B_k\frac{x^k}{k!}=1-\frac{x}2+\sum_{k=1}^\infty B_{2k}\frac{x^{2k}}{(2k)!}, \quad |x|<2\pi.
\end{equation}
In combinatorics, Stirling numbers of the second kind $S(n,k)$ for $n\ge k\ge0$ may be computed by
\begin{equation}\label{Stirling-Number-dfn}
S(n,k)=\frac1{k!}\sum_{\ell=0}^k(-1)^{k-\ell}\binom{k}{\ell}\ell^{n}
\end{equation}
and may be generated by
\begin{equation}\label{2stirling-gen-funct-exp}
\frac{(e^x-1)^k}{k!}=\sum_{n=k}^\infty S(n,k)\frac{x^n}{n!}, \quad k\in\{0\}\cup\mathbb{N}.
\end{equation}
\par
In~\cite[p.~536]{GKP-Concrete-Math-1989} and~\cite[p.~560]{GKP-Concrete-Math-2nd}, the following simple formula for computing Bernoulli numbers $B_n$ in terms of Stirling numbers of the second kind $S(n,k)$ was incidentally obtained.

\begin{theorem}\label{Bernoulli-Stirling-thm}
For $n\in\{0\}\cup\mathbb{N}$, we have
\begin{equation}\label{Bernoulli-Stirling-eq}
B_n=\sum_{k=0}^n(-1)^k\frac{k!}{k+1}S(n,k).
\end{equation}
\end{theorem}

The aim of this paper is to provide four alternative proofs for the explicit formula~\eqref{Bernoulli-Stirling-eq}.

\section{Four alternative proofs of the formula~\eqref{Bernoulli-Stirling-eq}}

Now we start out to provide four alternative proofs for the explicit formula~\eqref{Bernoulli-Stirling-eq}.
\par
Considering $S(0,0)=1$, it is clear that the formula~\eqref{Bernoulli-Stirling-eq} is valid for $n=0$.  Further considering $S(n,0)=0$ for $n\ge1$, it is sufficient to show
\begin{equation}\label{Bernoulli-Stirling-nozero}
B_n=\sum_{k=1}^n(-1)^k\frac{k!}{k+1}S(n,k), \quad n\in\mathbb{N}.
\end{equation}

\begin{proof}[First proof]
It is listed in~\cite[p.~230, 5.1.32]{abram} that
\begin{equation}\label{ln-frac}
\ln\frac{b}a=\int_0^\infty\frac{e^{-au}-e^{-bu}}u\td u.
\end{equation}
Taking $a=1$ and $b=1+x$ in~\eqref{ln-frac} yields
\begin{equation}\label{ln-frac-div-x-int}
\frac{\ln(1+x)}x=\int_0^\infty\frac{1-e^{-xu}}{xu}e^{-u}\td u
=\int_0^\infty\biggl(\int_{1/e}^1t^{xu-1}\td t\biggr)e^{-u}\td u.
\end{equation}
Replacing $x$ by $e^x-1$ in~\eqref{ln-frac-div-x-int} results in
\begin{equation}\label{Bernoulli-Gen-Funct-Int-Eq}
\frac{x}{e^x-1}=\int_0^\infty\biggl(\int_{1/e}^1t^{ue^x-u-1}\td t\biggr)e^{-u}\td u.
\end{equation}
\par
In combinatorics, Bell polynomials of the second kind, or say, the partial Bell polynomials, $\bell_{n,k}(x_1,x_2,\dotsc,x_{n-k+1})$ are defined by
\begin{equation}
\bell_{n,k}(x_1,x_2,\dotsc,x_{n-k+1})=\sum_{\substack{1\le i\le n,\ell_i\in\mathbb{N}\\ \sum_{i=1}^ni\ell_i=n\\ \sum_{i=1}^n\ell_i=k}}\frac{n!}{\prod_{i=1}^{n-k+1}\ell_i!} \prod_{i=1}^{n-k+1}\Bigl(\frac{x_i}{i!}\Bigr)^{\ell_i}
\end{equation}
for $n\ge k\ge1$, see~\cite[p.~134, Theorem~A]{Comtet-Combinatorics-74}, and satisfy
\begin{equation}\label{Bell(n-k)}
\bell_{n,k}\bigl(abx_1,ab^2x_2,\dotsc,ab^{n-k+1}x_{n-k+1}\bigr) =a^kb^n\bell_{n,k}(x_1,x_n,\dotsc,x_{n-k+1})
\end{equation}
and
\begin{equation}\label{Bell-stirling}
\bell_{n,k}(\overbrace{1,1,\dotsc,1}^{n-k+1})=S(n,k),
\end{equation}
see~\cite[p.~135]{Comtet-Combinatorics-74}, where $a$ and $b$ are any complex numbers. The well-known Fa\`a di Bruno formula may be described in terms of Bell polynomials of the second kind \linebreak $\bell_{n,k}(x_1,x_2,\dotsc,x_{n-k+1})$ by
\begin{equation}\label{Bruno-Bell-Polynomial}
\frac{\td^n}{\td x^n}f\circ g(x)=\sum_{k=1}^nf^{(k)}(g(x)) \bell_{n,k}\bigl(g'(x),g''(x),\dotsc,g^{(n-k+1)}(x)\bigr),
\end{equation}
see~\cite[p.~139, Theorem~C]{Comtet-Combinatorics-74}.
\par
Applying in~\eqref{Bruno-Bell-Polynomial} the functions $f(y)=t^{y}$ and $g(x)=ue^x-u-1$ gives
\begin{equation}\label{n-deriv-(uexp)}
\frac{\td^nt^{ue^x}}{\td x^n}=\sum_{k=1}^n(\ln t)^kt^{ue^x} \bell_{n,k}\bigl(\overbrace{ue^x,ue^x,\dotsc,ue^x}^{n-k+1}\bigr).
\end{equation}
Making use of the formulas~\eqref{Bell(n-k)} and~\eqref{Bell-stirling} in~\eqref{n-deriv-(uexp)} reveals
\begin{equation}\label{n-deriv-(uexp)-S}
\frac{\td^nt^{ue^x}}{\td x^n}=t^{ue^x}\sum_{k=1}^nS(n,k)u^k(\ln t)^k e^{kx}.
\end{equation}
Differentiating $n$ times on both sides of~\eqref{Bernoulli-Gen-Funct-Int-Eq} and considering~\eqref{n-deriv-(uexp)-S} figure out
\begin{equation}\label{bernoulli-n-deriv-S}
\frac{\td^n}{\td x^n}\biggl(\frac{x}{e^x-1}\biggr)
=\sum_{k=1}^nS(n,k)e^{kx}\int_0^\infty u^k \biggl(\int_{1/e}^1 (\ln t)^kt^{ue^x-u-1} \td t\biggr)e^{-u}\td u.
\end{equation}
On the other hand, differentiating $n$ times on both sides of~\eqref{Bernumber-dfn} gives
\begin{equation}\label{Bernumber-dfn-deriv}
\frac{\td^n}{\td x^n}\biggl(\frac{x}{e^x-1}\biggr)=\sum_{k=n}^\infty B_k\frac{x^{k-n}}{(k-n)!}.
\end{equation}
Equating~\eqref{bernoulli-n-deriv-S} and~\eqref{Bernumber-dfn-deriv} and taking the limit $x\to0$ discover
\begin{align*}
B_n&=\sum_{k=1}^nS(n,k)\int_0^\infty u^k \biggl(\int_{1/e}^1 \frac{(\ln t)^k}{t} \td t\biggr)e^{-u}\td u\\
&=\sum_{k=1}^n\frac{(-1)^k}{k+1}S(n,k)\int_0^\infty u^k e^{-u}\td u\\
&=\sum_{k=1}^n\frac{(-1)^kk!}{k+1}S(n,k).
\end{align*}
The first proof of Theorem~\ref{Bernoulli-Stirling-thm} is complete.
\end{proof}

\begin{proof}[Second proof]
In the book~\cite[p.~386]{bullenmean} and in the papers~\cite[p.~615]{Carlson-AMM-72} and~\cite[p.~885]{Neuman-jmaa.1994.1469}, it was given that
\begin{equation}\label{log-mean-int-eq}
\frac{\ln b-\ln a}{b-a}=\int_0^1\frac{1}{(1-t)a+tb}\td t,
\end{equation}
where $a,b>0$ and $a\ne b$. Replacing $a$ by $1$ and $b$ by $e^x$ yields
\begin{equation}\label{log-mean-int-exp}
\frac{x}{e^x-1}=\int_0^1\frac1{1+(e^x-1)t}\td t.
\end{equation}
Applying the functions $f(y)=\frac1y$ and $y=g(x)=1+(e^x-1)t$ in the formula~\eqref{Bruno-Bell-Polynomial} and simplifying by~\eqref{Bell(n-k)} and~\eqref{Bell-stirling} give
\begin{align*}
\frac{\td^n}{\td x^n}\biggl(\frac{x}{e^x-1}\biggr)&=\int_0^1\frac{\td^n}{\td x^n} \biggl[\frac1{1+(e^x-1)t}\biggr]\td t\\*
&=\int_0^1\sum_{k=1}^n(-1)^k\frac{k!}{[1+(e^x-1)t]^{k+1}} \bell_{n,k}(\overbrace{te^x,te^x,\dotsc,te^x}^{n-k+1})\td t\\
&=\sum_{k=1}^n(-1)^kk!\int_0^1\frac{t^k}{[1+(e^x-1)t]^{k+1}} \bell_{n,k}(\overbrace{e^x,e^x,\dotsc,e^x}^{n-k+1})\td t\\
&\to\sum_{k=1}^n(-1)^kk!\int_0^1t^k \bell_{n,k}(\overbrace{1,1,\dotsc,1}^{n-k+1})\td t,\quad x\to0\\
&=\sum_{k=1}^n(-1)^kk!S(n,k)\int_0^1t^k\td t\\
&=\sum_{k=1}^n(-1)^k\frac{k!}{k+1}S(n,k).
\end{align*}
On the other hand, taking the limit $x\to0$ in~\eqref{Bernumber-dfn-deriv} leads to
\begin{equation*}
\frac{\td^n}{\td x^n}\biggl(\frac{x}{e^x-1}\biggr)=\sum_{k=n}^\infty B_k\frac{x^{k-n}}{(k-n)!}
\to B_n,\quad x\to0.
\end{equation*}
The second proof of Theorem~\ref{Bernoulli-Stirling-thm} is thus complete.
\end{proof}

\begin{proof}[Third proof]
Let $CT[f(x)]$ be the coefficient of $x^0$ in $f(x)$. Then
\begin{align*}
\sum_{k=1}^n(-1)^k \frac{k!}{k+1}S(n,k)
&=\sum_{k=1}^n(-1)^k CT\biggl[\frac{n!}{x^n}\frac{(e^x-1)^k}{k+1}\biggl]\\
&=n! CT \biggl[ \frac1{x^n} \sum_{k=1}^\infty(-1)^k\frac{(e^x-1)^k}{k+1}\biggr]\\
&=n! CT \biggl[ \frac1{x^n}\frac{\ln[1+(e^x-1)]- (e^x-1)}{e^x-1}\biggl]\\
&=n! CT \biggl[ \frac1{x^n} \frac{x}{e^x-1}\biggl]\\*
&=B_n.
\end{align*}
Thus, the formula~\eqref{Bernoulli-Stirling-eq} follows.
\end{proof}

\begin{proof}[Fourth proof]
It is clear that the equation~\eqref{Bernumber-dfn} may be rewritten as
\begin{equation}\label{Bernumber-dfn-rew}
\frac{\ln[1+(e^x-1)]}{e^x-1}=\sum_{k=0}^\infty B_k\frac{x^k}{k!}.
\end{equation}
Differentiating $n$ times on both sides of~\eqref{Bernumber-dfn-rew} and taking the limit $x\to0$ reveal
\begin{align*}
B_n&=\lim_{x\to0}\sum_{k=n}^\infty B_k\frac{x^{k-n}}{(k-n)!}
=\lim_{x\to0}\frac{\td^n}{\td x^n}\biggl(\frac{\ln[1+(e^x-1)]}{e^x-1}\biggr)\\
&=\lim_{x\to0}\sum_{k=1}^n\biggl[\frac{\ln(1+u)}{u}\biggr]^{(k)} \bell_{n,k}(\overbrace{e^x,e^x,\dotsc,e^x}^{n-k+1}), \quad u=e^x-1\\
&=\lim_{x\to0}\sum_{k=1}^n\Biggl[\sum_{\ell=1}^\infty(-1)^{\ell-1}\frac{u^{\ell-1}}{\ell}\Biggr]^{(k)} \bell_{n,k}(\overbrace{e^x,e^x,\dotsc,e^x}^{n-k+1})\\
&=\lim_{x\to0}\sum_{k=1}^n\Biggl[\sum_{\ell=k+1}^\infty(-1)^{\ell-1} \frac{(\ell-1)!}{(\ell-k-1)!\ell}u^{\ell-k-1}\Biggr] \bell_{n,k}(\overbrace{e^x,e^x,\dotsc,e^x}^{n-k+1})\\
&=\sum_{k=1}^n\lim_{u\to0} \Biggl[\sum_{\ell=k+1}^\infty(-1)^{\ell-1} \frac{(\ell-1)!}{(\ell-k-1)!\ell}u^{\ell-k-1}\Biggr] \lim_{x\to0}\bell_{n,k}(\overbrace{e^x,e^x,\dotsc,e^x}^{n-k+1})\\
&=\sum_{k=1}^n(-1)^k\frac{k!}{k+1} \bell_{n,k}(\overbrace{1,1,\dotsc,1}^{n-k+1})\\
&=\sum_{k=1}^n(-1)^k\frac{k!}{k+1}S(n,k).
\end{align*}
The fourth proof of Theorem~\ref{Bernoulli-Stirling-thm} is thus complete.
\end{proof}

\begin{remark}
In~\cite[p.~1128, Corollary]{recursion}, among other things, it was found that
\begin{equation}\label{Bernoulli-N-Guo-Qi-99}
B_{2k}= \frac12 - \frac1{2k+1} - 2k \sum_{i=1}^{k-1}
\frac{A_{2(k-i)}}{2(k - i) + 1}
\end{equation}
for $k\in\mathbb{N}$, where $A_m$ is defined by
\begin{equation*}
\sum_{m=1}^nm^k=\sum_{m=0}^{k+1}A_mn^{m}.
\end{equation*}
\par
In~\cite[p.~559]{GKP-Concrete-Math-2nd} and~\cite[Theorem~2.1]{exp-derivative-sum-Combined.tex}, it was collected and recovered that
\begin{equation}\label{exp-deriv-exp}
\biggl(\frac1{e^x-1}\biggr)^{(k)}=(-1)^k\sum_{m=1}^{k+1}(m-1)!S(k+1,m)\biggl(\frac1{e^x-1}\biggr)^m, \quad k\in\{0\}\cup\mathbb{N}.
\end{equation}
In~\cite[Theorem~3.1]{exp-derivative-sum-Combined.tex}, by the identity~\eqref{exp-deriv-exp}, it was obtained that
\begin{multline}\label{Bernumber-formula-eq}
B_{2k}=1+\sum_{m=1}^{2k-1}\frac{S(2k+1,m+1) S(2k,2k-m)}{\binom{2k}{m}} \\*
-\frac{2k}{2k+1}\sum_{m=1}^{2k}\frac{S(2k,m)S(2k+1,2k-m+1)}{\binom{2k}{m-1}}, \quad k\in\mathbb{N}.
\end{multline}
\par
In~\cite[Theorem~1.4]{Tan-Cot-Bernulli-No.tex}, among other things, it was presented for $k\in\mathbb{N}$ that
\begin{equation}\label{Bernoulli-N-Explicit}
B_{2k}=\frac{(-1)^{k-1}k}{2^{2(k-1)}(2^{2k}-1)}\sum_{i=0}^{k-1}\sum_{\ell=0}^{k-i-1} (-1)^{i+\ell}\binom{2k}{\ell}(k-i-\ell)^{2k-1}.
\end{equation}
\par
Recently, a new formula
\begin{equation}\label{Bernoulli-Stirling-formula}
B_n=\sum_{i=0}^n(-1)^{i}\frac{\binom{n+1}{i+1}}{\binom{n+i}{i}}S(n+i,i)
\end{equation}
for $n\in\mathbb{N}$ was discovered in the preprint~\cite{Bernoulli-No-Int-New.tex}.
\end{remark}

\begin{remark}
The identities in~\eqref{exp-deriv-exp} have been generalized and applied in~\cite{Eight-Identy-More.tex, CAM-D-13-01430-Xu-Cen}.
\end{remark}

\begin{remark}
This paper is a slightly revised version of the preprint~\cite{Bernoulli-Stirling2-3P.tex}.
\end{remark}

\subsection*{Acknowledgements}
The authors thank Professor Doron Zeilberg in USA for his reminding of the books~\cite{GKP-Concrete-Math-1989, GKP-Concrete-Math-2nd} and sketching the third proof in an e-mail on October 10, 2013. Due to this, the authors find that the formula~\eqref{Bernoulli-Stirling-eq} originated from the uneasily-found literature~\cite{Logan-Bell-87} and was listed as an incidental consequence of an answer to an exercise in~\cite[p.~536]{GKP-Concrete-Math-1989} and~\cite[p.~560]{GKP-Concrete-Math-2nd}.

\end{document}